\documentclass[amstex,12pt]{article}
\usepackage{amsfonts}
\usepackage{amssymb}
\usepackage{amsmath}
\usepackage{mathrsfs}
\usepackage{amssymb,amsmath,amsthm}
\newcommand{\beq}{\begin{equation}}
\newcommand{\eeq}{\end{equation}}
\newcommand{\beqn}{\begin{eqnarray}}
\newcommand{\eeqn}{\end{eqnarray}}

\def\e{\hat{\mathbb{E}}}
\def\h{\mathcal {H}}

\def\o{{\omega}}
\renewcommand\O{{\Omega}}

\begin{document}
\small
\title{Stability theorems for stochastic differential equations driven by  G-Brownian motion}
\author{{Defei Zhang}\thanks{Corresponding author. E-mail: zhdefei@163.com.}
\\ \small{ Department of Mathematics,
Shandong University, Jinan  250100, China}\\
 }
\date{}
\maketitle
\begin{abstract} In this paper,  stability theorems for stochastic differential equations
 and backward stochastic differential equations driven by G-Brownian motion
are  obtained.  We  show the  existence  and  uniqueness of   solutions to forward-backward stochastic differential equations driven by  G-Brownian motion. Stability theorem for forward-backward stochastic differential equations driven by  G-Brownian motion is also presented.
\end{abstract}
\section{Introduction}
Consider a family of ordinary stochastic differential equations (SDEs for short) parameterized by $\varepsilon\geq 0,$
$$X_t^\varepsilon=x_0^\varepsilon+\int_0^tb^\varepsilon(s,X_s^\varepsilon)ds+\int_0^t\sigma^\varepsilon(s,X_s^\varepsilon)dW_s, t\in [0,T],$$
where $W_t$ is classical Brownian motion. It is well known that the
strong convergence of the coefficient in $L^2$ implies the strong
convergence of the solutions, that is, if \begin{center}
$x_0^\varepsilon\rightarrow x^0_0$, as $\varepsilon\rightarrow 0,$
\end{center}
and
\begin{center}
$E[\int_0^T(|b^\varepsilon(s,X_s^0)-b^0(s,X_s^0)|^2+|\sigma^\varepsilon(s,X_s^0)-\sigma^0(s,X_s^0)|^2)ds]\rightarrow
0,$ as $\varepsilon\rightarrow 0,$
\end{center}
then under Lipschitz and other reasonable assumptions, their
solutions also converge strongly in $L^2,$
\begin{center} $\forall \ t\in [0,T], E[|X_t^\varepsilon-X_t^0|^2]\rightarrow
0,$ as $\varepsilon\rightarrow 0.$
\end{center}
This result, known as the continuous dependence theorem, or the
stability property, can be found in many standard textbooks of SDEs
(e.g., see \cite{p6}).

 Backward stochastic differential equations
(BSDEs for short) driven by classical Brownian motion were introduced, in linear case, by Bismut
\cite{b1} in 1973.
  In 1990, Pardoux and Peng  considered general BSDEs (see \cite{p1}). Similar continuous dependence theorem
for the case of backward stochastic differential equations was
obtained by El Karoui,  Peng and Quenez (1994) \cite{e} and  Hu and Peng (1997)
\cite{h1}.

  As  for  the  forward-backward  equations, Antonelli  \cite{a} first  studied  these
equations,  and  he  gave  the  existence  and  uniqueness  when  the  time  duration  T
is  sufficiently  small. Using  a  PDE  approach,  Ma, Protter and Yong   \cite{m}  gave  the  existence
and  uniqueness  to  a  class  of  forward-backward  SDEs  in  which  the  forward
SDE  is  non-degenerate. In 1995, Hu and Peng \cite{h2} study  the  existence  and  uniqueness  of  the  solutions  to  forward-backward  stochastic  differential  equations  without  the  non-degeneracy condition.

Motivated by  uncertainty problems, risk measures and the
superhedging in finance, Peng(2006, see \cite{p2}) has  introduced  the notion of sublinear expectation space,
which is a generalization of classical probability space.
 Together with the notion of sublinear expectation, Peng also introduced the related G-normal distribution and G-Brownian motion.
 The expectation
 associated with G-Brownian motion is a sublinear expectation which
 is called G-expectation. The stochastic calculus with respect to the G-Brownian
 motion has been established by Peng in \cite{p2}, \cite{p3} and \cite{p4}. Since these notions were introduced, many properties of G-Brownian motion have been studied by authors, for example,  \cite{d}, \cite{g1}, \cite{g2}, \cite{s1}-\cite{x},  et al.

Therefore,  the natural  questions are: stability properties for stochastic differential equations and backward stochastic differential equations
driven by G-Brownian motion are also true?  How to obtain the existence  and  uniqueness of  the  solution of  a
 forward-backward stochastic differential equations driven by  G-Brownian motion?
The   goal  of this paper is to study stability properties for
stochastic differential equations driven by G-Brownian motion
(G-SDEs for short) and backward stochastic differential equations
driven by G-Brownian motion (G-BSDEs for short). Indeed, under
Lipschitz or integral-Lipschitz condition and  other reasonable assumptions,
stability theorems for G-SDEs and G-BSDEs are obtained. Meanwhile, we also show   the  existence  and  uniqueness of  the  solution of  a
new type of forward-backward stochastic differential equations driven by  G-Brownian motion.

 This paper is organized as follows: in Section 2, we recall briefly  some notions and properties
about G-expectation and G-Brownian motion. In Section 3, we study
the stability properties of G-SDEs, while  Section 4, study the
G-BSDEs case.  At last,  the  existence  and  uniqueness of  the  solution of   forward-backward stochastic differential equations driven by  G-Brownian motion are obtained.  Stability theorem for forward-backward stochastic differential equations driven by  G-Brownian motion is also presented.
\section{Preliminaries}
In this section, we introduce some notations and preliminaries about sublinear expectations and G-Brownian motion, which
will be needed in what follows. More details concerning this section may be found in \cite{p2}, \cite{p3} and \cite{p4}.

 Let $\Omega$ be a
given set and let $\mathcal {H}$ be a linear space of real valued
bounded functions defined on $\Omega$. We suppose that $\mathcal
{H}$ satisfies  $C\in\mathcal {H}$ for each constant
$C$ and $|X|\in\mathcal {H},$ if  $X\in\mathcal {H}.$\\
\textbf{Definition 2.1}. A sublinear expectation $\mathbb{E}$ is a
functional $\mathbb{E}:\mathcal {H}\rightarrow R$ satisfying
\\(i) Monotonicity: $\mathbb{E}[X]\geq \mathbb{E}[Y]$ if $X\geq
Y.$\\
(ii) Constant preserving: $\mathbb{E}[C]=C$ for $C\in R.$\\
(iii) Sub-additivity: For each
$X,Y\in\h,\mathbb{E}[X+Y]\leq\mathbb{E}[X]+\mathbb{E}[Y].$
\\(iv) Positive homogeneity: $\mathbb{E}[\lambda X]=\lambda\mathbb{E}[X]$ for $\lambda\geq
0.$

 The triple $(\O,\h,\mathbb{E})$ is called a sublinear
expectation space. If (i) and (ii) are satisfied,
$\mathbb{E}[\cdot]$ is called a nonlinear expectation and the triple
$(\O,\h,\mathbb{E})$ is called a nonlinear expectation space.

From now on, we consider the following sublinear expectation space
$(\Omega,{\cal H},\mathbb{E})$: if $X_1,\cdots,X_n\in{\cal H}$,
then $\varphi(X_1,\cdots,X_n)\in{\cal H}$ for each $\varphi\in
C_{l,Lip}(\mathbb R^n)$, where $C_{l,Lip}(R^{n})$ denotes the linear space of functions $\varphi$ satisfying $|\varphi(x)-\varphi(y)|\leq C(1+|x|^m+|y|^m)|x-y|$ for $x,y\in R^n$, some $C>0,m\in \mathbb{N}$ depending on $\varphi.$
\\ \textbf{Definition 2.2 }. Let $X$ and $Y$ be two $n$-dimensional random vectors defined on nonlinear expectation spaces  $(\O_1,\h_1,\mathbb{E}_1)$ and $(\O_2,\h_2,\mathbb{E}_2)$, respectively. They are called identically distributed, denoted by $X\stackrel{d}{=}Y$, if \begin{center}$\mathbb{E}_1[\varphi(X)]=\mathbb{E}_2[\varphi(Y)]$, for $\forall \varphi\in C_{l,Lip}(R^n).$ \end{center}
 \textbf{Definition 2.3}. In a nonlinear expectation space $(\O,\h,\mathbb{E})$, a random vector $Y\in \h^n$ is said to be independent from another random vector $X\in \h^m$ under $\mathbb{E}[\cdot]$, if
\begin{center}$\mathbb{E}[\varphi(X,Y)]=\mathbb{E}[\mathbb{E}[\varphi(x,Y)]_{x=X}]$, for $\forall \varphi\in C_{l,Lip}(R^{m+n}).$\end{center}
$\bar{X}$ is called an independent copy of $X$ if $\bar{X}\stackrel{d}{=}X$ and $\bar{X}$ is independent from $X.$\\
\textbf{Definition 2.4 (G-normal distribution)}. In a sublinear expectation space $(\O,\h,\mathbb{E})$, a random variable $X \in \h$ with
$$\mathbb{E}[X^2]=\bar{\sigma}^2, -\mathbb{E}[-X^2]=\underline{\sigma}^2,$$ is said to be $N(0; [\underline{\sigma}^2 , \bar{\sigma}^2 ])$-distributed, if for
each $\bar{X} \in \h$ which is an independent copy of  $X$ we have $$aX+b\bar{X}\stackrel{d}{=}\sqrt{a^2+b^2}X,\ \ \forall a,b\geq 0.$$
\textbf{Definition 2.5 (G-Brownian motion)}. A process $\{B_t(\o)\}_{t\geq 0}$ in a sublinear expectation space $(\O,\h,\mathbb{E})$,
is called a G-Brownian motion if for each $n \in N$ and $0 \leq t_1\leq \cdots \leq t_n <\infty,B_{t_1},\cdots ,B_{t_n}\in \h$ and the following properties are satisfied:\\
(i) $B_0(\o)=0;$\\
(ii) For each $t, s \geq 0,$ the increment $B_{t+s}- B_t$ is $N(0; [\underline{\sigma}^2 , \bar{\sigma}^2 ])$-distributed and is independent from $(B_{t_1},\cdots ,B_{t_n})$ for each $n \in N$ and $0 \leq t_1\leq\cdots \leq t_n \leq t.$

 We denote by ${{\Omega}}=C_0^d(R^+)$ the space
of all $R^d$-valued continuous paths $(\o_t)_{t\in R^+}$, with
$\o_0=0,$ equipped with the distance
$\rho(\o^1,\o^2):=\sum\limits_{i=1}^\infty
2^{-i}[(\max\limits_{t\in[0,i]}|\o_t^1-\o_t^2|)\wedge 1].$
Considering the canonical process $B_t(\o)=(\o_t)_{t\geq 0}$. For
each fixed $T>0$, set $\O_T:=\{\omega_{.\wedge T}:\omega\in \O\}$
and
$$L_{ip}(\O_T):=\{\varphi(B_{t_1},B_{t_2},...,B_{t_m}):m\geq 1,t_1,...,t_m\in[0,T],\varphi\in C_{l,Lip}(R^{d\times m})\},$$
and define $L_{ip} (\O):=\bigcup\limits_{n=1}^\infty L_{ip}(\O_n).$

 Let $\xi$ be a G-normal distributed, or  $N(0; [\sigma^2, 1])$-distributed random variable in a sublinear expectation space $(\widetilde{\O},\widetilde{\h},\widetilde{\mathbb{E}}).$ We now introduce a sublinear expectation $\e$ defined on $L_{ip} (\O)$ via the following procedure: for each $X\in L_{ip}(\O)$ with $$X=\varphi(B_{t_1}-B_{t_0},B_{t_2}-B_{t_1},\cdots,B_{t_m}-B_{t_{m-1}}),$$ for some $\varphi\in C_{l,Lip}(R^{d\times m})$ and $0=t_0<t_1<\cdots<t_m<\infty$, we set
$$\e[\varphi(B_{t_1}-B_{t_0},\cdots,B_{t_m}-B_{t_{m-1}})]:=\widetilde{\mathbb{E}}[\varphi(\sqrt{t_1-t_0}\xi_1,\cdots,\sqrt{t_m-t_{m-1}}\xi_m)],$$
where $(\xi_1,\cdots,\xi_m)$ is an m-dimensional G-normal distributed random vector
in a sublinear expectation space $(\widetilde{\O},\widetilde{\h},\widetilde{\mathbb{E}})$ such that $\xi_i\stackrel{d}{=} N(0; [\sigma^2, 1])$ and
such that $\xi_{i+1}$ is independent from $(\xi_1,\cdots,\xi_i)$ for each $i = 1,\cdots,m $.

The related conditional expectation of $X=\varphi(B_{t_1}-B_{t_0},B_{t_2}-B_{t_1},...,B_{t_m}-B_{t_{m-1}})$ under $\O_{t_j}$ is defined by
\begin{equation*}\begin{aligned}\e[X|\O_{t_j}]&=\e[\varphi(B_{t_1},B_{t_2}-B_{t_1},\cdots,B_{t_m}-B_{t_{m-1}})|\O_{t_j}]\\
&:=\psi(B_{t_1},,\cdots,B_{t_j}-B_{t_{j-1}}),
\end{aligned}\end{equation*}
where $$\psi(x_1,\cdots,x_j)=\widetilde{\mathbb{E}}[\varphi(x_1,\cdots,x_j,\sqrt{t_{j+1}-t_j}\xi_{j+1},\cdots,\sqrt{t_m-t_{m-1}}\xi_m)].$$
\textbf{Definition 2.6}. The expectation $\e[\cdot]:L_{ip}(\O)\rightarrow R$
 defined through the
above procedure is called G-expectation. The corresponding canonical process $(B_t)_{t\geq 0}$ in the sublinear expectation space $(\O,L_{ip}(\O),\e)$ is called a G-Brownian motion.

 We denote by $L_G^p(\O_T),p\geq 1,$ the completion of
$L_{ip}(\O_T) $ under the norm $||X||_p:=(\e[|X|^p])^{1/p}.$
Similarly, denote $L_G^p(\O)$ is complete space of $L_{ip}(\O).$ We
give some important properties about conditional G-expectation
$\e[\cdot|\O_t],t\in[0,T].$\\
 \textbf{Proposition 2.1}. The
conditional expectation $\e[\cdot|\O_t],t\in[0,T]$ holds for each
$X,Y\in L_G^1(\O_t):$\\
(i) If $X\geq Y,$ then $\e[X|\O_t]\geq \e[Y|\O_t].$\\
(ii) $\e[\eta|\O_t]=\eta,$ for each $t\in[0,\infty)$ and $\eta\in
L_G^1(\O_t).$
\\(iii) $\e[X|\O_t]-\e[Y|\O_t]\leq \e[X-Y|\O_t].$\\
(iv)
 $\e[\eta X|\O_t]=\eta^+\e[X|\O_t]+\eta^-\e[-X|\O_t]$ for each bounded
 $\eta\in L_G^1(\O_t). $\\
 (v) $\e[\e[X|\O_t]|\O_s]=\e[X|\O_{t\wedge s}]$, in particular,
 $\e[\e[X|\O_t]]=\e[X].$

Next, we introduce the It\^{o}'s integral with G-Brownian motion. For
$T\in R^+,$ a partition ${\pi}_T$ of $[0,T]$ is a finite ordered
subset $\pi_T=\{t_0,t_1,...,t_N\}$ such that $0=t_0<t_1<...<t_N=T,$
$$\mu(\pi_T):=\max\{|t_{i+1}-t_i|:i=0,1,...,N-1\}.$$

 Using
$\pi_T^N=\{t_0^N,t_1^N,...,t_N^N\}$ to denote a sequence of
partitions of $[0,T]$ such that
$\lim\limits_{N\rightarrow\infty}\mu(\pi_T^N)=0.$\\
Let $p\geq 1$
 be fixed. We consider the following type of simple processes: for a given
 partition $\pi_T=\{t_0,t_1,...,t_N\}$ of $[0,T],$ set $\eta_t(\o)=\sum\limits_{k=0}^{N-1}\xi_k(\o)\mathbf{I}_{[t_k,t_{k+1})}(t),$
 where $\xi_k\in L_G^P(\O_{t_k}),k=0,1,...,N-1$ are given. The collection of these processes
 is denoted by $M_G^{p,0}(0,T).$ For each $p\geq 1$, we denote by  $M_G^p([0,T]; R^n)$  the completion of
$M_G^{p,0}([0,T]; R^n)$ under the norm
$||\eta_t||_{M_G^p([0,T])}:=(\int_0^T\e[|\eta_t|^p]dt)^{\frac{1}{p}}.$\\
\textbf{Definition 2.7}. For an $\eta\in M_G^{p,0}(0,T),$ the related
Bochner integral is
$$\int_0^T\eta_t(\o)dt:=\sum\limits_{k=0}^{N-1}\xi_k(\o)(t_{k+1}-t_k).$$

Let $(B_t)_{t\geq 0}$ be a 1-dimensional G-Brownian motion with
$G(a):=\frac{1}{2}\e[aB_1^2]=\frac{1}{2}(\bar{\sigma}^2a^+-\underline{\sigma}^2a^-),$ where $\bar{\sigma}^2=\e[B_1^2], \underline{\sigma}^2=-\e[-B_1^2],$
 $0\leq \underline{\sigma}\leq\bar{\sigma}<\infty.$\\
\textbf{Definition 2.8}. For an $\eta\in M_G^{2,0}(0,T)$ of the form
$\eta_t(\o)=\sum\limits_{k=0}^{N-1}\xi_k(\o)\mathbf{I}_{[t_k,t_{k+1})}(t),$
define
$$\int_0^T\eta(s)dB_s:=\sum\limits_{k=0}^{N-1}\xi_k(B_{t_{k+1}}-B_{t_k}).$$
\textbf{Proposition 2.2}. For each $\eta\in M_G^{2,0}(0,T),$ then
$$\e[\int_0^T\eta(s)dB_s]=0,\ \
\e[(\int_0^T\eta(s)dB_s)^2]\leq\bar{{\sigma}}^2\int_0^T\e[\eta^2(s)]ds.$$
\textbf{Definition 2.9}. For the 1-dimensional G-Brownian motion
$B_t$, we denote $\langle B\rangle_t$ is the quadratic variation
process of $B_t,$ where $\langle
B\rangle_t:=\lim\limits_{\mu(\pi_t^N)\rightarrow
0}\sum\limits_{k=0}^{N-1}(B_{t_{k+1}^N}-B_{t_k^N})^2=(B_{t})^2-2\int_0^tB_{s}dB_{s}.$
\\
\textbf{Definition 2.10}. For each $\eta\in M_G^{1,0}(0,T),$  define
$\int_0^T\eta(s)d\langle
B\rangle_s:=\sum\limits_{k=0}^{N-1}\xi_k(\langle
B\rangle_{t_{k+1}}-\langle B\rangle_{t_k}).$
\\ \textbf{Proposition 2.3}. For any $0\leq t\leq
T<\infty,$\\
(i) $\e[|\int_0^T\eta_td\langle B\rangle_t|]\leq
\bar{{\sigma}}^2\e[\int_0^T|\eta_t|dt],\forall\  \eta_t\in
M_G^{1}(0,T).$ \\(ii)
$\e[(\int_0^T\eta_tdB_t)^2]=\e[\int_0^T\eta_t^2d\langle
B\rangle_t],\forall \ \eta_t\in M_G^{2}(0,T).$\\
(iii) $\e[\int_0^T|\eta_t|^pdt]\leq \int_0^T\e[|\eta_t|^p]dt, \forall
\ {\eta}_t\in M_G^{p}(0,T),\ p\geq 1.$
\section{Stability theorem of G-stochastic differential equations}
 In this
section, we consider the
stability theorem  of G-stochastic differential equations. Consider
the following stochastic differential equations driven by
d-dimensional G-Brownian motion:
\begin{equation} X_t=X_0+\int_0^t
b(s,X_s)ds+\sum\limits_{i,j=1}^d\int_0^t h_{ij}(s,X_s)d\langle
B^i,B^j\rangle_s+\sum\limits_{j=1}^d\int_0^t\sigma_j(s,X_s)dB_s^j,\
\ t\in[0,T],
\end{equation}
the initial condition $X_0\in R^n$, and $b,
h_{ij},\sigma_j$ are given functions satisfying $b(\cdot,x),
h_{ij}(\cdot,x),$ $\sigma_j(\cdot, x)\in M_G^2([0,T]; R^n)$ for each
$x\in R^n.$  Consider the
following G-SDEs depending on a parameter
$\varepsilon(\varepsilon\geq 0)$:
\begin{equation} X^\varepsilon_t=X^\varepsilon_0+\int_0^t
b^\varepsilon(s,X_s^\varepsilon)ds+\sum\limits_{i,j=1}^d\int_0^t
h^\varepsilon_{ij}(s,X_s^\varepsilon)d\langle
B^i,B^j\rangle_s+\sum\limits_{j=1}^d\int_0^t\sigma^\varepsilon_j(s,X^\varepsilon_s)dB_s^j,\
\ t\in[0,T].
\end{equation}

We make the following assumptions:\\
 \textbf{Assumption 3.1}. For any $\varepsilon \geq 0, x \in R^n, b^\varepsilon(\cdot, x), h_{ij}^\varepsilon(\cdot,
 x),
 \sigma_j^\varepsilon(\cdot, x)\in M_G^2([0,T]; R^n), X_0^\varepsilon\in
 R^n.$\\
\textbf{Assumption 3.2}. For any $\varepsilon \geq 0, x, x_1, x_2\in R^n:$\\
(H1)
$|b^\varepsilon(t,x)|^2+\sum\limits_{i,j=1}^d|h^\varepsilon_{ij}(t,x)|^2+\sum\limits_{j=1}^d|\sigma^\varepsilon_{j}(t,x)|^2\leq
\alpha_1^2(t)+\alpha_2^2(t)|x|^2,$\\
(H2) $|b^\varepsilon(t,x_1)-b^\varepsilon(t,x_2)|^2+\sum\limits_{i,j=1}^d|h^\varepsilon_{ij}(t,x_1)-h^\varepsilon_{ij}(t,x_2)|^2+\sum\limits_{j=1}^d|\sigma^\varepsilon_{j}(t,x_1)-\sigma^\varepsilon_{j}(t,x_2)|^2\leq
\alpha^2(t)\rho(|x_1-x_2|^2),$ where $\alpha_1\in M_G^2([0,T]),
 \alpha_2:[0,T]\rightarrow R^+$ and $\alpha:[0,T]\rightarrow R^+$ are
Lebesgue integrable, and $\rho:(0,+\infty)\rightarrow(0,+\infty)$ is
continuous, increasing, concave function satisfying $\rho(0+)=0,
\int_0^1\frac{1}{\rho(r)}dr=+\infty.$\\
\textbf{Assumption 3.3}. (i) $\forall \ t\in [0,T],$ as
$\varepsilon\rightarrow 0,$
$$\int_0^t\e[|\phi^\varepsilon(s,X^0_s)-\phi^0(s,X^0_s)|^2]ds\rightarrow
0,$$where $\phi=b, h_{ij}$ and $\sigma_j,$ respectively, $i,j=1,\cdots,d$.\\
 (ii) As $\varepsilon\rightarrow 0,$
$$X_0^\varepsilon\rightarrow X_0^0.$$
\textbf{Remark 3.1}. The Assumptions 3.1 and 3.2 guarantee, for any
$\varepsilon\geq 0,$ the existence of a unique solution
$X_t^\varepsilon\in M_G^2([0,T];R^n)$ of G-SDEs (3.2)(see
\cite{b}), while the Assumption 3.3 will allow us to deduce the
following stability
theorem for G-SDEs.\\
 \textbf{Theorem 3.1}. Under the Assumptions 3.1, 3.2 and 3.3, we have the
 following convergence: as $\varepsilon\rightarrow 0,$
\begin{equation}
\forall \ t\in[0,T],\ \ \e[|X_t^\varepsilon-X_t^0|^2]\rightarrow 0.
\end{equation}

In order to prove Theorem 3.1, we need the following lemmas:\\
\textbf{Lemma 3.1} (see Chemin  and Lerner \cite{c}). Let
$\rho:(0,+\infty)\rightarrow(0,+\infty)$ be a continuous, increasing
function satisfying  $\rho(0+)=0,$ $
\int_0^1\frac{1}{\rho(r)}dr=+\infty $ and let $u$ be a measurable,
nonnegative function defined on $(0,+\infty)$ satisfying
$$u(t)\leq a+\int_0^t\beta(s)\rho(u(s))ds,\ t\in (0,+\infty),$$
where $a\in[0,+\infty),$  and $\beta:[0,T]\rightarrow R^+$ is
Lebesgue integrable. Then\\
(i) if $a=0$, then $u(t)=0$, for $t\in[0,+\infty);$\\
(ii) if $a>0$, then $$u(t)\leq v^{-1}(v(a)+\int_0^t\beta(s)ds),$$
where $v(t):=\int_{t_0}^t\frac{1}{\rho(s)}ds, t_0\in (0,+\infty).$\\
\textbf{Lemma 3.2} (see Peng \cite{p4}). Let $\rho: R\rightarrow R$ be a
continuous increasing, concave function defined on $R,$ then for
each  $X\in L_G^1(\O), \forall \ t\geq 0,$ the following
Jensen inequality holds:
$$\rho(\e[X|\O_t])\geq \e[\rho(X)|\O_t].$$
 \textbf{Proof of Theorem 3.1}. Let $\hat{X}_t^\varepsilon:=X_t^\varepsilon-X_t^0,$
 $\hat{X}_0^\varepsilon:=X_0^\varepsilon-X_0^0,$ then
\begin{equation}\begin{aligned}
\hat{X}_t^\varepsilon&=\hat{X}_0^\varepsilon+\int_0^t
(b^\varepsilon(s,X_s^\varepsilon)-b^0(s,X_s^0))ds\\&+\sum\limits_{i,j=1}^d\int_0^t
(h^\varepsilon_{ij}(s,X_s^\varepsilon)-h^0_{ij}(s,X_s^0))d\langle
B^i,B^j\rangle_s\\&+\sum\limits_{j=1}^d\int_0^t(\sigma^\varepsilon_j(s,X^\varepsilon_s)-\sigma^0_j(s,X^0_s))dB_s^j,
\end{aligned}\end{equation} and
\begin{equation}\begin{aligned}
|\hat{X}_t^\varepsilon|^2&\leq
C\{|\hat{X}_0^\varepsilon|^2+|\int_0^t
(b^\varepsilon(s,X_s^\varepsilon)-b^\varepsilon(s,X_s^0))ds|^2+|\int_0^t
(b^\varepsilon(s,X_s^0)-b^0(s,X_s^0))ds|^2\\&+\sum\limits_{i,j=1}^d|\int_0^t
(h^\varepsilon_{ij}(s,X_s^\varepsilon)-h^\varepsilon_{ij}(s,X_s^0))d\langle
B^i,B^j\rangle_s|^2+\sum\limits_{i,j=1}^d|\int_0^t
(h^\varepsilon_{ij}(s,X_s^0)-h^0_{ij}(s,X_s^0))d\langle
B^i,B^j\rangle_s|^2\\&+\sum\limits_{j=1}^d|\int_0^t(\sigma^\varepsilon_j(s,X^\varepsilon_s)-\sigma^\varepsilon_j(s,X^0_s))dB_s^j|^2
+\sum\limits_{j=1}^d|\int_0^t(\sigma^\varepsilon_j(s,X^0_s)-\sigma^0_j(s,X^0_s))dB_s^j|^2\},
\end{aligned}\end{equation}
taking the G-expectation on both sides of the above relation and from
Proposition 2.3, we get
\begin{equation}\begin{aligned}
\e[|\hat{X}_t^\varepsilon|^2]&\leq
C\{|\hat{X}_0^\varepsilon|^2+\int_0^t\e[|
b^\varepsilon(s,X_s^\varepsilon)-b^\varepsilon(s,X_s^0)|^2]ds+\int_0^t\e[|
b^\varepsilon(s,X_s^0)-b^0(s,X_s^0)|^2]ds\\&+\sum\limits_{i,j=1}^d\int_0^t\e[|
h^\varepsilon_{ij}(s,X_s^\varepsilon)-h^\varepsilon_{ij}(s,X_s^0)|^2]ds+\sum\limits_{i,j=1}^d\int_0^t\e[|
h^\varepsilon_{ij}(s,X_s^0)-h^0_{ij}(s,X_s^0)|^2]ds\\&+\sum\limits_{j=1}^d\int_0^t\e[|\sigma^\varepsilon_j(s,X^\varepsilon_s)-\sigma^\varepsilon_j(s,X^0_s)|^2]ds
+\sum\limits_{j=1}^d\int_0^t\e[|\sigma^\varepsilon_j(s,X^0_s)-\sigma^0_j(s,X^0_s)|^2]ds\},
\end{aligned}\end{equation}
by  Assumption 3.2, we  have
\begin{equation}\begin{aligned}
\e[|\hat{X}_t^\varepsilon|^2]\leq
C^\varepsilon(T)+C_2\int_0^t\alpha^2(s)\e[\rho(|\hat{X}_s^\varepsilon|^2)]ds,
\end{aligned}\end{equation}
where
\begin{equation*}\begin{aligned}C^\varepsilon(t):&=C\int_0^t\e[|
b^\varepsilon(s,X_s^0)-b^0(s,X_s^0)|^2]ds+
C\sum\limits_{i,j=1}^d\int_0^t\e[|
h^\varepsilon_{ij}(s,X_s^0)-h^0_{ij}(s,X_s^0)|^2]ds\\&+C\sum\limits_{j=1}^d
\int_0^t\e[|\sigma^\varepsilon_j(s,X^0_s)-\sigma^0_j(s,X^0_s)|^2]ds+C|\hat{X}_0^\varepsilon|^2.\end{aligned}\end{equation*}
Because $\rho$ is concave and increasing,  from Lemma 3.2, we have
\begin{equation}\begin{aligned} \e[|\hat{X}_t^\varepsilon|^2]\leq
C^\varepsilon(T)+C_2\int_0^t\alpha^2(s)\rho(\e[|\hat{X}_s^\varepsilon|^2])ds.
\end{aligned}\end{equation}
Since as
$\varepsilon\rightarrow 0,$ $C^\varepsilon(T)\rightarrow 0,$  hence,
from  Lemma 3.1, we get
\begin{center}
$ \e[|\hat{X}_t^\varepsilon|^2]\rightarrow 0,$ as
$\varepsilon\rightarrow 0.$\end{center}The proof is
complete.

A special case of  Assumption 3.2 is\\
\textbf{Assumption 3.4 (Lipschitz condition)}. For any $ x_1, x_2 \in R^n,$ there exist
constant $C_0>0$ such that
$$|\phi^\varepsilon(t,x_1)-\phi^\varepsilon(t,x_2)|\leq C_0 |x_1-x_2|, \  \ t\in[0,T],$$ where
$\phi=b, h_{ij}$ and $\sigma_j,$ respectively, $i,j=1,\cdots,d$. \\
 \textbf{Corollary 3.1}. Under the Assumptions 3.1, 3.3 and 3.4, we have the
 convergence of the solution of the G-SDEs (3.2) in the sense
of (3.3).
\section{Stability theorem of G-backward stochastic differential equations}
 \setcounter{equation}{0}
 In this section, we give a stability theorem of backward stochastic
differential equations  driven by  d-dimensional G-Brownian
motion (G-BSDEs for short). Consider the following type of G-backward
stochastic differential equations depending on a parameter
$(\delta\geq 0)$:
\begin{equation}\begin{aligned}
Y_t^\delta=\e[\xi^\delta+\int_t^Tf^\delta(s,Y^\delta_s)ds+\sum\limits_{i,j=1}^d\int_t^Tg^\delta_{ij}(s,Y^\delta_s)d\langle
B^i, B^j\rangle_s|\O_t],\ t\in [0,T],
\end{aligned}\end{equation}
where $\xi^\delta\in L_G^1(\O_T; R^n)$ is given, and  $
f^\delta(\cdot,y), g^\delta_{ij}(\cdot,y)\in
M_G^1(0,T; R^n).$

We further make the following assumptions:
\\
\textbf{Assumption 4.1}. For any $\delta\geq 0,y, y_1, y_2 \in
R^n,$\\(H1)
$|f^\delta(t,y)|+\sum\limits_{i,j=1}^d|g^\delta_{ij}(t,y)|\leq
\beta(t)+C|y|,$\\
(H2)$|f^\delta(t,y_1)-f^\delta(t,y_2)|+\sum\limits_{i,j=1}^d|g^\delta_{ij}(t,y_1)-g^\delta_{ij}(t,y_2)|\leq
\rho(|y_1-y_2|),$\\ where $C>0, \beta\in M_G^1([0,T];R^+),$  and
$\rho:(0,+\infty)\rightarrow(0,+\infty)$ is continuous, increasing,
concave function satisfying $\rho(0+)=0,
\int_0^1\frac{1}{\rho(r)}dr=+\infty.$\\
\textbf{Assumption 4.2}. (i) $\forall \ t\in [0,T],$ as
$\delta\rightarrow 0,$
$$\int_t^T\e[|\phi^\delta(s,Y^0_s)-\phi^0(s,Y^0_s)|]ds\rightarrow
0,$$ where $\phi=f, g_{ij}$  respectively, $i,j=1,\cdots,d$.\\
 (ii) As $\delta\rightarrow 0,$
$$\e[|\xi^\delta-\xi^0|]\rightarrow 0.$$
\textbf{Remark 4.1}. Under the Assumptions 4.1  and 4.2,  G-BSDEs (4.1) has a unique solution. The proof goes in a similar way as that in \cite{b}, and we omit it.
\\ \textbf{Theorem 4.1}. Under the Assumptions 4.1 and 4.2, we have the
 following convergence: as $\delta\rightarrow 0,$
\begin{equation}
\forall \ t\in[0,T],\ \ \e[|Y_t^\delta-Y_t^0|]\rightarrow 0.
\end{equation}
\textbf{Proof.} Let $\hat{Y}_t^\delta:=Y_t^\delta-Y_t^0,
\hat{\xi}^\delta:=\xi^\delta-\xi^0,$ then
\begin{equation}\label{eq:4.14}\begin{aligned}
|\hat{Y}^\delta_t|&\leq
\e[|\hat{\xi}^\delta|+\int_t^T|f^\delta(s,Y^\delta_s)-f^0(s,Y^0_s)|ds+\sum\limits_{i,j=1}^d\int_t^T|g^\delta_{ij}(s,Y^\delta_s)-g^0_{ij}(s,Y^0_s)|d\langle
B^i, B^j\rangle_s|\O_t]\\
&\leq
\e[|\hat{\xi}^\delta|+\int_t^T|f^\delta(s,Y^0_s)-f^0(s,Y^0_s)|ds+\sum\limits_{i,j=1}^d\int_t^T|g^\delta_{ij}(s,Y^0_s)-g^0_{ij}(s,Y^0_s)|d\langle
B^i, B^j\rangle_s\\
&+\int_t^T|f^\delta(s,Y^\delta_s)-f^\delta(s,Y^0_s)|ds+\sum\limits_{i,j=1}^d\int_t^T|g^\delta_{ij}(s,Y^\delta_s)-g^\delta_{ij}(s,Y^0_s)|d\langle
B^i, B^j\rangle_s|\O_t].\\
\end{aligned}\end{equation}
Taking the G-expectation on both sides of (\ref{eq:4.14}), we have
\begin{equation}\begin{aligned}
\e[|\hat{Y}^\delta_t|]&\leq
\e[|\hat{\xi}^\delta|]+\int_t^T\e[|f^\delta(s,Y^\delta_s)-f^0(s,Y^0_s)|]ds+C\sum\limits_{i,j=1}^d\int_t^T\e[|g^\delta_{ij}(s,Y^\delta_s)-g^0_{ij}(s,Y^0_s)|]ds\\
&+\int_t^T\e[|f^\delta(s,Y^\delta_s)-f^\delta(s,Y^0_s)|]ds+C\sum\limits_{i,j=1}^d\int_t^T\e[|g^\delta_{ij}(s,Y^\delta_s)-g^\delta_{ij}(s,Y^0_s)|]ds.
\end{aligned}\end{equation}
From the Assumption 4.1, Propositions 2.1 and 2.3  as well as Lemma 3.2, we have
\begin{equation}\begin{aligned}
\e[|\hat{Y}^\delta_t|]&\leq C^\delta(0)+
K_1\int_t^T\e[\rho(|\hat{Y}^\delta_s|)]ds
\\&\leq C^\delta(0)+
K_1\int_t^T\rho(\e[|\hat{Y}^\delta_s|])ds.
\end{aligned}\end{equation}
where $$
C^\delta(0):=\e[|\hat{\xi}^\delta|]+\int_0^T\e[|f^\delta(s,Y^0_s)-f^0(s,Y^0_s)|]ds+C\sum\limits_{i,j=1}^d\int_0^T\e[|g^\delta_{ij}(s,Y^0_s)-g^0_{ij}(s,Y^0_s)|]ds.$$
Since as $\delta\rightarrow 0,$ $C^\delta(0)\rightarrow 0,$ hence, from Lemma 3.1, we have
$$ \e[|\hat{Y}^\delta_t|]\rightarrow 0.$$ The proof is
complete.

A special case of  Assumption 4.1 is\\
\textbf{Assumption 4.3}. For any $\delta\geq 0, y_1, y_2 \in R^n,$
there exist constant $C_0>0$ such that
$$|\phi^\delta(t,y_1)-\phi^\delta(t,y_2)|\leq C_0 |y_1-y_2|,\  t\in[0,T],$$
$\phi=f, g_{ij}$  respectively, $i,j=1,\cdots,d$. \\
 \textbf{Corollary 4.1}. Under the Assumptions 4.2  and 4.3, we have the
 convergence of the solution of the G-BSDEs (4.1) in the sense
of (4.2).
\section{Forward-backward stochastic differential equations}
 \setcounter{equation}{0}
The  goal  of  this  section is  to  show  the  existence  and
uniqueness of forward-backward  stochastic differential equations driven by G-Brownian motion.
 For notational simplification, we only consider the case of 1-dimensional G-Brownian motion. However,  our method can be easily extend to the case of multi-dimensional G-Brownian motion. We consider the following system:
\begin{equation} \label{eq:3.1}
\left\{ \begin{aligned}
 X_t&=x+\int_0^t
b(s,X_s,Y_s)ds+\int_0^t h(s,X_s,Y_s)d\langle
B\rangle_s+\int_0^t\sigma(s,X_s,Y_s)dB_s,\\
Y_t&=\e[\xi+\int_t^Tf(s,X_s,Y_s)ds+\int_t^Tg(s,X_s,Y_s)d\langle
B\rangle_s|\O_t],\ t\in [0,T],
                          \end{aligned} \right.
                          \end{equation}
where the initial condition $x\in R$, the terminal data $\xi\in
L_G^2(\O_T;R)$, and $b, h,\sigma, f,g $ are given functions
satisfying $b(\cdot,x,y), h(\cdot,x,y),$ $\sigma(\cdot,
x,y),f(\cdot,x,y),g(\cdot,x,y)\in M_G^2([0,T]; R)$ for any
$(x,y)\in R^2$ and the Lipschitz condition, i.e.,
$|\phi(t,x,y)-\phi(t,x',y')|\leq K(|x-x'|+|y-y'|),$ for each
$t\in[0,T],(x,y)\in R^2, (x',y')\in R^2,\phi=b, h,\sigma, f$ and
$g$, respectively. The solution is a pair of processes $(X,Y)\in
M_G^2(0,T;R)\times M_G^2(0,T;R).$

 This model is  called forward-backward because the  two components
  in  the  system (\ref{eq:3.1}) are  solutions,  respectively,  of  a  G-forward and
  a G-backward stochastic differential equation.

We first introduce the following mappings on a fixed interval
$[0,T]:$
$$\Lambda^i_{\cdot}:M_G^2(0,T;R)\times M_G^2(0,T;R)\rightarrow M_G^2(0,T;R)\times M_G^2(0,T;R),i=1,2,$$
by setting $\Lambda^i_{t},i=1,2, t\in[0,T],$ with
\begin{equation} \label{eq:3.2}
 \begin{aligned}
\Lambda^1_{t}(X,Y)&=x+\int_0^t b(s,X_s,Y_s)ds+\int_0^t
h(s,X_s,Y_s)d\langle
B\rangle_s+\int_0^t\sigma(s,X_s,Y_s)dB_s,\\
\Lambda^2_{t}(X,Y)&=\e[\xi+\int_t^Tf(s,X_s,Y_s)ds+\int_t^Tg(s,X_s,Y_s)d\langle
B\rangle_s|\O_t].
                          \end{aligned}
                          \end{equation}
 \textbf{Lemma 5.1}. For any $(X,Y),(X',Y')\in
M_G^2(0,T;R)\times M_G^2(0,T;R),$ we have the following estimates:
\begin{equation} \label{eq:3.3}
 \begin{aligned}
\e[|\Lambda^1_{t}(X,Y)-\Lambda^1_{t}(X',Y')|^2]&\leq C \int_0^t
\e[|X_s-X'_s|^2+|Y_s-Y'_s|^2]ds,t\in [0,T],\\
\e[|\Lambda^2_{t}(X,Y)-\Lambda^2_{t}(X',Y')|^2]&\leq C' \int_t^T
\e[|X_s-X'_s|^2+|Y_s-Y'_s|^2]ds,t\in [0,T],
                          \end{aligned}
                          \end{equation}
where $C=24K^2, C'=8K^2, K$ is  Lipschitz  coefficient.\\
 Proof.  \begin{equation*}\begin{aligned}
\e[|\Lambda^1_{t}(X,Y)-\Lambda^1_{t}(X',Y')|^2]&\leq
4\int_0^t\e[|b(s,X_s,Y_s)-b(s,X'_s,Y'_s)|^2]ds\\
&+4\int_0^t\e[|h(s,X_s,Y_s)-h(s,X'_s,Y'_s)|^2]ds\\
&+4\int_0^t\e[|\sigma(s,X_s,Y_s)-\sigma(s,X'_s,Y'_s)|^2]ds\\
&\leq 24K^2\int_0^t \e[|X_s-X'_s|^2+|Y_s-Y'_s|^2]ds.
                          \end{aligned}
                          \end{equation*}
And since \begin{equation*}
 \begin{aligned}
|\Lambda^2_{t}(X,Y)-\Lambda^2_{t}(X',Y')|^2&\leq
2\e[|\int_t^Tf(s,X_s,Y_s)-f(s,X'_s,Y'_s)ds|^2\\
&+|\int_t^Tg(s,X_s,Y_s)-g(s,X'_s,Y'_s)ds|^2|\O_t],
                          \end{aligned}
                          \end{equation*}
then
 \begin{equation*}
 \begin{aligned}
\e[|\Lambda^2_{t}(X,Y)-\Lambda^2_{t}(X',Y')|^2]&\leq
2\e[|\int_t^Tf(s,X_s,Y_s)-f(s,X'_s,Y'_s)ds|^2\\
&+|\int_t^Tg(s,X_s,Y_s)-g(s,X'_s,Y'_s)ds|^2]\\
&\leq 2 \int_t^T\e[|f(s,X_s,Y_s)-f(s,X'_s,Y'_s)|^2]ds\\
&+2\int_t^T\e[|g(s,X_s,Y_s)-g(s,X'_s,Y'_s)|^2]ds\\
&\leq 8K^2 \int_t^T \e[|X_s-X'_s|^2+|Y_s-Y'_s|^2]ds.
                          \end{aligned}
                          \end{equation*}

 Let  us  consider  the  space $M_G^2(0,T;R)\times M_G^2(0,T;R),$  with  the  norm
                           $||(X,Y)||_{M_G^2(0,T)\times
                           M_G^2(0,T)}:=||X||_{M_G^2(0,T)}+||Y||_{M_G^2(0,T)}=\int_0^T\e[|X_s|^2]ds+\int_0^T\e[|Y_s|^2]ds,$
                           this is a Banach space.
\\ \textbf{Theorem  5.1}. Under  the  previously  stated  hypotheses  and  under  the additional  assumption
that $(2\sqrt{6}+2\sqrt{2})K\sqrt{T}<1$, then there exists a unique
solution $(X,Y)\in M_G^2(0,T;R)\times M_G^2(0,T;R)$ of the
forward-backward stochastic differential equation (\ref{eq:3.1}) in
the
${M_G^2(0,T)\times M_G^2(0,T)}$ sense. \\
Proof. Let  us  consider  the  space $M_G^2(0,T;R)\times
M_G^2(0,T;R),$ with  the  norm $$||(X,Y)||_{M_G^2(0,T)\times
M_G^2(0,T)}:=||X||_{M_G^2(0,T)}+||Y||_{M_G^2(0,T)}.$$  We can view
the  system (\ref{eq:3.1})  as the operator $\Lambda_t(X,Y):= \left(
  \begin{array}{c}
    \Lambda^1_t(X,Y)\\
      \Lambda^2_t(X,Y)\\
  \end{array}
\right),$
 thus
\begin{equation}\label{eq:3.6}
 \begin{aligned}
&||\Lambda_t(X,Y)-\Lambda_t(X',Y')||_{M_G^2(0,T)\times M_G^2(0,T)}\\&=
||\Lambda^1_t(X,Y)-\Lambda^1_t(X',Y')||_{M_G^2(0,T)}+||\Lambda^2_t(X,Y)-\Lambda^2_t(X',Y')||_{M_G^2(0,T)}\\
&=(\int_0^T
\e[(\Lambda^1_t(X,Y)-\Lambda^1_t(X',Y'))^2]dt)^{\frac{1}{2}}\\&+(\int_0^T
\e[(\Lambda^2_t(X,Y)-\Lambda^2_t(X',Y'))^2]dt)^{\frac{1}{2}}.
                          \end{aligned}
                          \end{equation}
Because the terminal data $\xi\in L_G^2(\O_T;R)$, and $b, h,\sigma,
f,g $ are given functions satisfying $b(\cdot,x,y), h(\cdot,x,y),$
$\sigma(\cdot, x,y),f(\cdot,x,y),g(\cdot,x,y)\in M_G^2([0,T]; R)$
for any $(x,y)\in R^2$ and the Lipschitz condition, we can prove
$||\Lambda_t(X,Y)||_{M_G^2(0,T)\times M_G^2(0,T)}<+\infty, \forall
(X,Y)\in{M_G^2(0,T;R)\times M_G^2(0,T;R)}.$ Next, we prove it is a
contraction mapping.
 From the Lemma 5.1, we can obtain
\begin{equation} \label{eq:3.7}
 \begin{aligned}
&||\Lambda_t(X,Y)-\Lambda_t(X',Y')||_{M_G^2(0,T)\times M_G^2(0,T)}\\
&\leq(\int_0^TC \int_0^t
\e[|X_s-X'_s|^2+|Y_s-Y'_s|^2]dsdt)^\frac{1}{2}\\&+(\int_0^T C' \int_t^T\e[|X_s-X'_s|^2+|Y_s-Y'_s|^2]dsdt)^\frac{1}{2}\\
&\leq
(\sqrt{C}+\sqrt{C'})\sqrt{T}\int_0^T\e[|X_s-X'_s|^2+|Y_s-Y'_s|^2]ds\\&=(\sqrt{C}+\sqrt{C'})\sqrt{T}||(X-X',Y-Y')||_{M_G^2(0,T)\times
M_G^2(0,T)}.
                          \end{aligned}
                          \end{equation}
From the assumption $(2\sqrt{6}+2\sqrt{2})K\sqrt{T}<1$, we can
obtain that $\Lambda_t(X,Y)$ is a contraction mapping. Hence a
unique fixed point for $\Lambda$ exists and this is  the solution of
our system (\ref{eq:3.1}).  The proof is complete.

In the last section, we present stability  theorem  for  forward-backward stochastic differential equations driven by G-Brownian motion.

\section{Stability  theorem  of  forward-backward stochastic differential equations}
 \setcounter{equation}{0}
Consider a family of forward-backward stochastic differential equations with parameter $(\gamma\geq 0),$
\begin{equation} \label{eq:6.1}
\left\{ \begin{aligned}
 X_t^\gamma&=x^\gamma+\int_0^t
b^\gamma(s,X^\gamma_s,Y^\gamma_s)ds+\int_0^t h^\gamma(s,X^\gamma_s,Y^\gamma_s)d\langle
B\rangle_s+\int_0^t\sigma^\gamma(s,X^\gamma_s,Y^\gamma_s)dB_s,\\
Y^\gamma_t&=\e[\xi^\gamma+\int_t^Tf^\gamma(s,X^\gamma_s,Y^\gamma_s)ds+\int_t^Tg^\gamma(s,X^\gamma_s,Y^\gamma_s)d\langle
B\rangle_s|\O_t],\ t\in [0,T],
                          \end{aligned} \right.
                          \end{equation}
where the initial condition $x^\gamma\in R$, the terminal data $\xi^\gamma\in
L_G^2(\O_T;R)$, and $b^\gamma, h^\gamma,\sigma^\gamma, f^\gamma,g^\gamma $ are given functions
satisfying $b^\gamma(\cdot,x,y), h^\gamma(\cdot,x,y),$ $\sigma^\gamma(\cdot,
x,y),f^\gamma(\cdot,x,y),g^\gamma(\cdot,x,y)\in M_G^2([0,T]; R)$ for any
$(x,y)\in R^2$ and the Lipschitz condition, i.e.,\\
\textbf{Assumption 6.1}.
$$|\phi(t,x,y)-\phi(t,x',y')|\leq K(|x-x'|+|y-y'|),$$ for each
$t\in[0,T],(x,y)\in R^2, (x',y')\in R^2,\phi=b^\gamma, h^\gamma,\sigma^\gamma, f^\gamma$ and
$g^\gamma$, respectively.

We further make the following assumption:\\
\textbf{Assumption 6.2}. (i) $\forall \ t\in [0,T],$ as
$\gamma\rightarrow 0,$
$$\int_0^t\e[|\phi^\gamma(s,X^0_s, Y_s^0)-\phi^0(s,X^0_s, Y_s^0)|^2]ds\rightarrow
0,$$where $\phi=b, h, \sigma, f$ and $g$, respectively.\\
 (ii) As $\gamma\rightarrow 0,$
$$x^\gamma\rightarrow x^0,\ \  \e[|\xi^\gamma-\xi^0|^2]\rightarrow 0.$$
\textbf{Theorem 6.1}. Under the Assumptions 6.1  and 6.2, then as $\gamma\rightarrow 0,$ $(X_t^\gamma, Y_t^\gamma)$ convergence to  $(X_t^0, Y_t^0)$ in the sense that
\begin{equation}
\forall \ t\in[0,T],\ \ \e[|X_t^\gamma-X_t^0|^2+|Y_t^\gamma-Y_t^0|^2]\rightarrow 0.
\end{equation}

The proof of Theorem 6.1 is similar to that of the Theorem 3.1; so, we omit it.
\section*{Acknowledgements}
The author thanks Professor Zengjing Chen for helpful discussion and suggestion.

\end{document}